# An aperiodic tiling of variable geometry made of two tiles, a triangle and a rhombus of any angle.

Vincent Van Dongen, PhD

October 17, 2021

Vincent.vandongen@gmail.com

## Abstract
Aperiodic tiling is a well-know area of research. First developed by mathematicians for the mathematical challenge they represent and the beauty of their resulting patterns, they became a growing field of interest when their practical use started to emerge. This was mainly in the eighties when a link was established with quasi-periodic materials. Several aperiodic tilings made of two tiles were discovered, the first one being by Penrose in the seventies. Since then, scientists discovered other aperiodic tilings including the square-triangle one, a tiling that has been particularly useful for the study of dodecagonal quasicrystals and soft matters. Based on this previous work, we discovered an infinite number of aperiodic tilings made of two tiles, a triangle and a rhombus of any angle. As a result, a variable geometry, i.e. continuously transformable, aperiodic tiling is proposed, whose underlying structure is dodecagonal. We discuss this limit case where the rhombus is so thin that it becomes invisible. At the boundary of this infinite space of tilings are two periodic ones; this represents a uniform view of periodic and aperiodic tilings.

## About this paper
This paper presents a new aperiodic tiling made of two tiles: a triangle and a rhombus of any angle. The figure below shows one of its instance when the rhombus has a 10-degree angle.

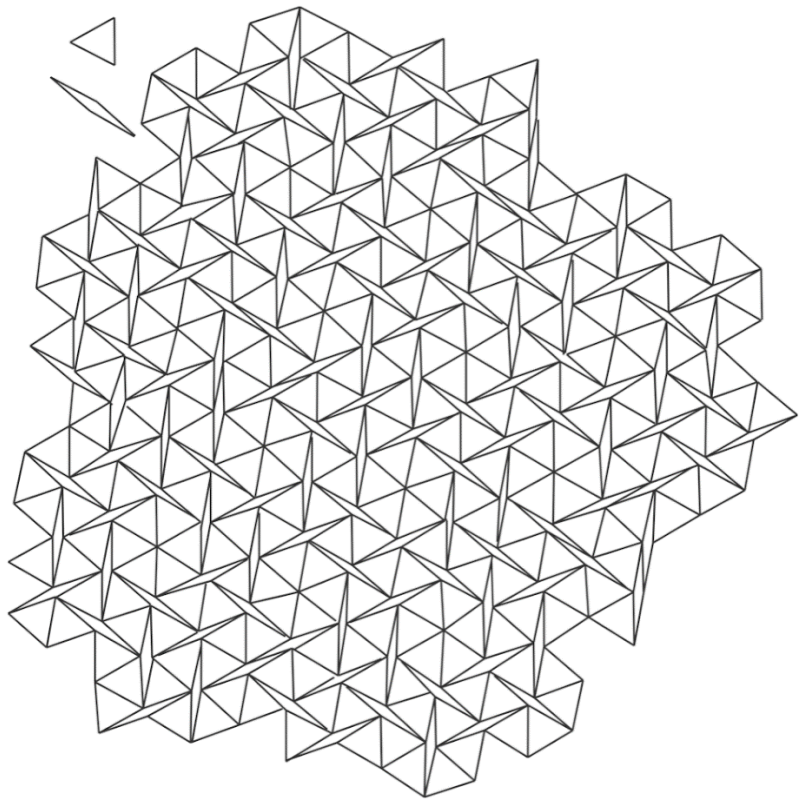

*Figure 1: An aperiodic tiling made of 2 tiles: a triangle and a rhombus of any angle – here a 10-degree angle.*

The presentation of this new tiling is organized as follows. In the next section, we recall the concepts of aperiodic tiling and some well-known examples. In particular, the square-triangle tiling (Baake, Klitzing, & Schlottmann, 1992) is presented and its link with dodecagonal quasicrystals. In the next section, we show how to create an infinite number of aperiodic tilings made of two tiles, a triangle and a rhombus (of any angle) whose edge is of the same length as the one of the triangle. The square-triangle tiling is a particular case in this family, where the rhombus is a square. We also show that in this infinite set of aperiodic tilings, at its boundaries, 2 of them are periodic. Finally, we discuss the case when the rhombus is so thin that it becomes invisible.

At the end, we conclude this paper with some perspectives for the future.

## Background information on aperiodic tiling

A lot of literature exists on aperiodic tiling and their various applications. In fact, whole books are covering this subject, as for example the one entitled *'Aperiodic Order'* (Baake & Grimm, 2013) (Baake & Grimm, Aperiodic Order - Vol. 2 Crystallography and Almost Periodicity, 2017). Here is a definition taken from MathWorld (Weisstein, s.d.): *An aperiodic tiling is a non-periodic tiling in which arbitrarily large periodic patches do not occur. A set of tiles is said to be aperiodic if they can form only non-periodic tilings. The most widely known examples of aperiodic tilings are those formed by Penrose tiles.*

Aperiodic tiling became popular when their link with quasi-periodic materials was established in the eighties. At this time, people like Prof. Dan Shechtman, who received the Nobel Prize of Chemistry in 2011, had seen the existence of quasi-periodic materials through TEM, transmission electron microscopy (Shechtman, s.d.). A link was then established with the first 2-tile aperiodic tiling discovered by Penrose that exhibits a fivefold symmetry (Penrose, 1979)

In 2004, already more than 8000 different quasicrystals were found (Steurer, 2004). Today, aperiodic structures are continuously being used and studied in various fields including, to name a few: aperiodic metal-organic frameworks (Oppenheim, Skorupskii, & Dincă, 2020), information theory (Varn & Crutchfield, 2016) and antenna placement optimization (Jian & Yang, 2020).

One may say that aperiodic tiling is easy to achieve. One can take an equilateral triangle for example and change its orientation in a random manner. The result will then be aperiodic. But we are interested instead by aperiodic tilings with structured patterns. One way of developing aperiodic tilings is by using substitution rules. Let us for example recall the square-triangle aperiodic tiling (Baake, Klitzing, & Schlottmann, 1992). A possible representation of this tiling is shown here below. (Ref: [Tilings Encyclopedia | Square-triangle (uni-bielefeld.de)](Tilings Encyclopedia | Square-triangle (uni-bielefeld.de)).

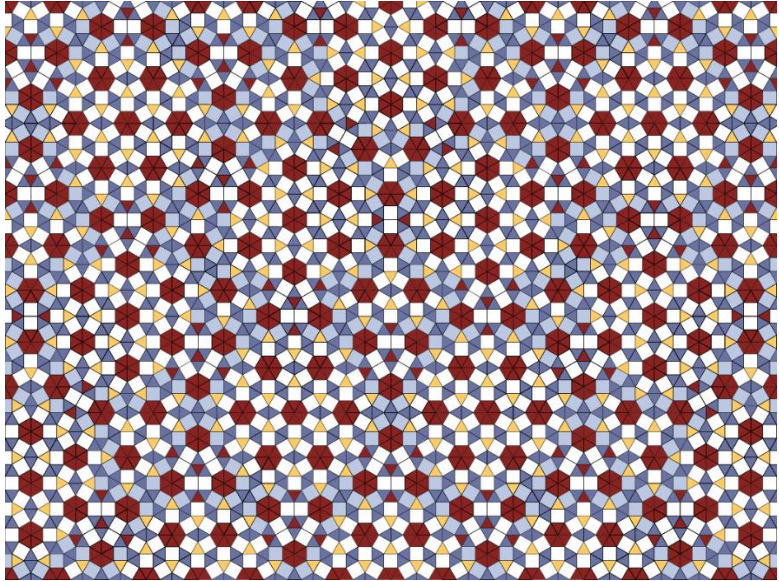

*Figure 2: Square-triangle aperiodic (Ref: Tilings Encyclopedia | Square-triangle (uni-bielefeld.de).*

The aperiodic tiling is generated by using 5 substitution rules. (Note that one of the tile, the light-blue square, is also used in flip mode.) The tiles and matching rules are shown here.

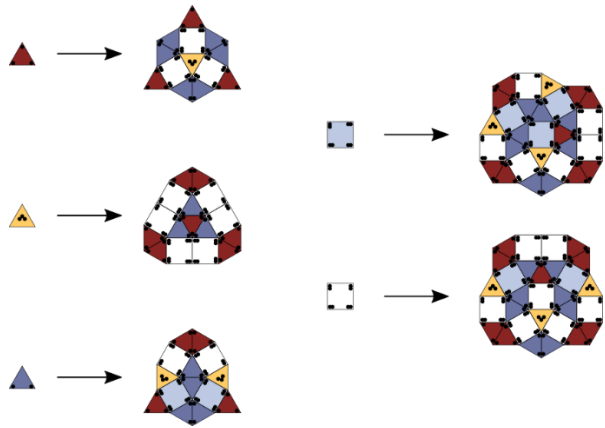

*Figure 3: Substitution rules for the square-triangle tiling (Ref: Tilings Encyclopedia | Square-triangle (uni-bielefeld.de).*

The matching rules on the tiles are needed to make use of the substitution rule and force the tiling to be aperiodic (and therefore structured). But once generated, these colors and patterns can be removed and a possible result is shown below. Hence, this aperiodic tiling makes use of two different tiles: a square and an equilateral triangle.

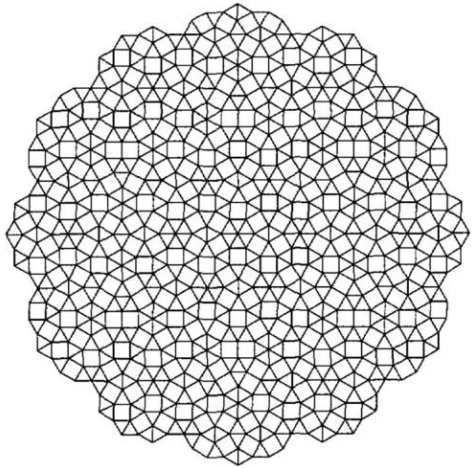

*Figure 4: The square-triangle aperiodic tiling is found in soft-matter quasicrystals (Image taken from (Sadoc & Mosseri))*

It has been shown that the square-triangle aperiodic tiling exhibits a 12-fold symmetry that can be found in non-metallic quasicrystals including soft-matter quasicrystals (Sadoc & Mosseri). These can also be referred as dodecagonal quasicrystals; the dodecagon can indeed be constructed by assembling equilateral triangles and squares. The square-triangle aperiodic tiling was recently generalized by using the same tiles but with new ways of assembling them (Impéror-Clerc, Jagannathan, Kalugin, & Sadoc, 2021).

## Introduction to variable-geometry aperiodic tiling

First, we observe that the square-triangle aperiodic tiling can be transformed into an infinite number of aperiodic tilings made of two tiles: a triangle (that is equilateral) and a rhombus. Indeed, if one forces one of the squares to change its angles, all other squares can follow the same change and the result is a different aperiodic tiling defined by these two tiles: a triangle and a rhombus. When the rhombus is a square, all angles of the rhombus are equal to 90 degrees and the tiling is the aperiodic square-triangle tiling. But the angle of the rhombus can take any value in the interval ]0, 180[. There is a special case, when the angle is 0, the rhombus disappears and the tiling becomes periodic, made of a single tile: an equilateral triangle. The same periodicity occurs when the angle is 180 degrees. This is illustrated here below.

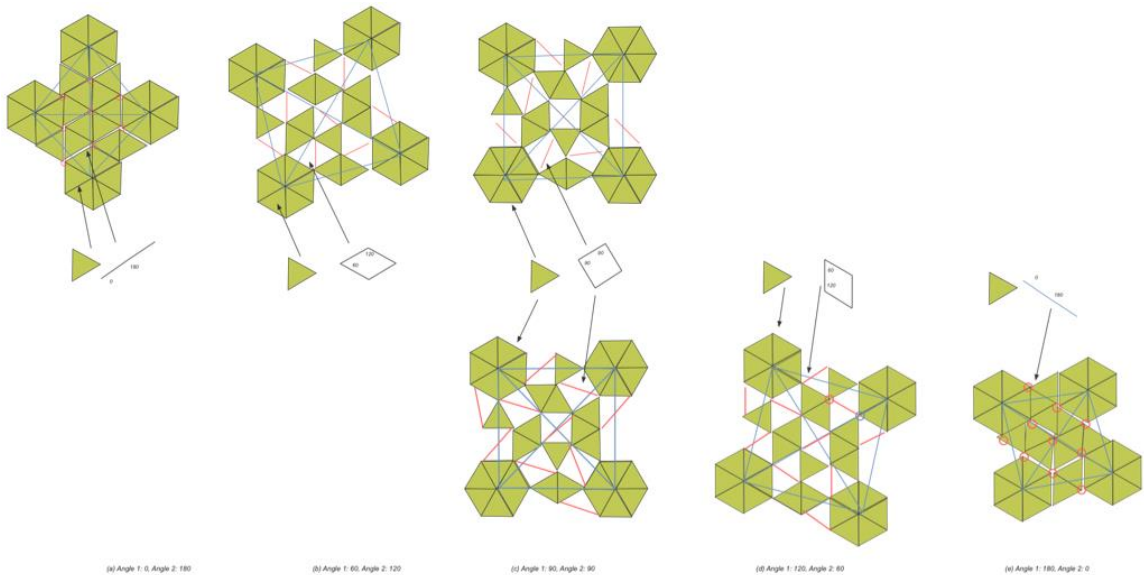

*Figure 5: Our new proposed family of tilings based on a triangle and a rhombus of any angle.*

All rhombi, for a given tiling, will always have the same angle because of the constraints imposed by the groups of triangles that are tight together by the vertices. Here below are the 5 possible sets of triangles based on the substitution rules.

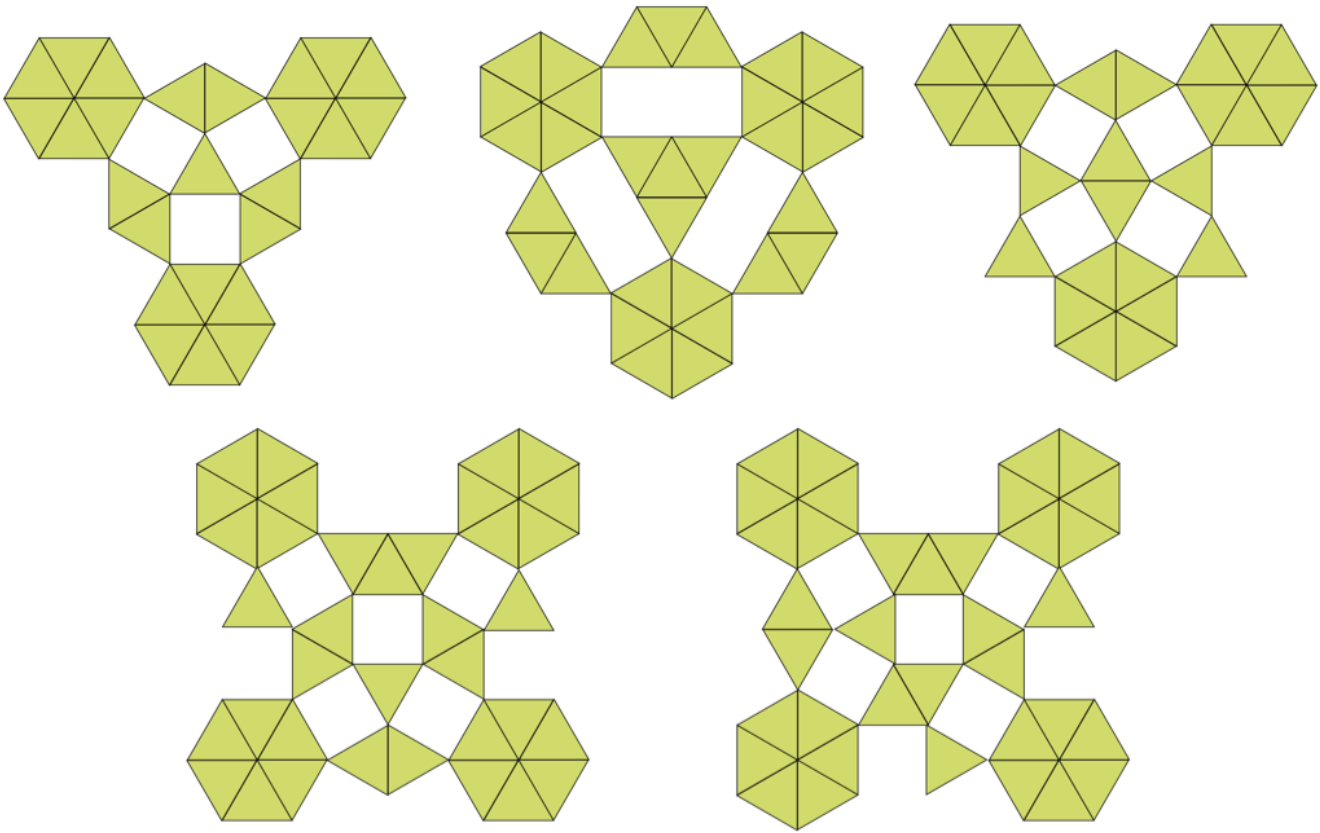

*Figure 6: Associated to the substitution rules of the square-triangle tiling, these combinations of triangle tiles attached by their vertices can act like pivots to allow for variable geometry (these combinations can also be flipped).*

These 5 combinations can then be used to created a large set of triangles attached by their vertices. An example is shown here below.

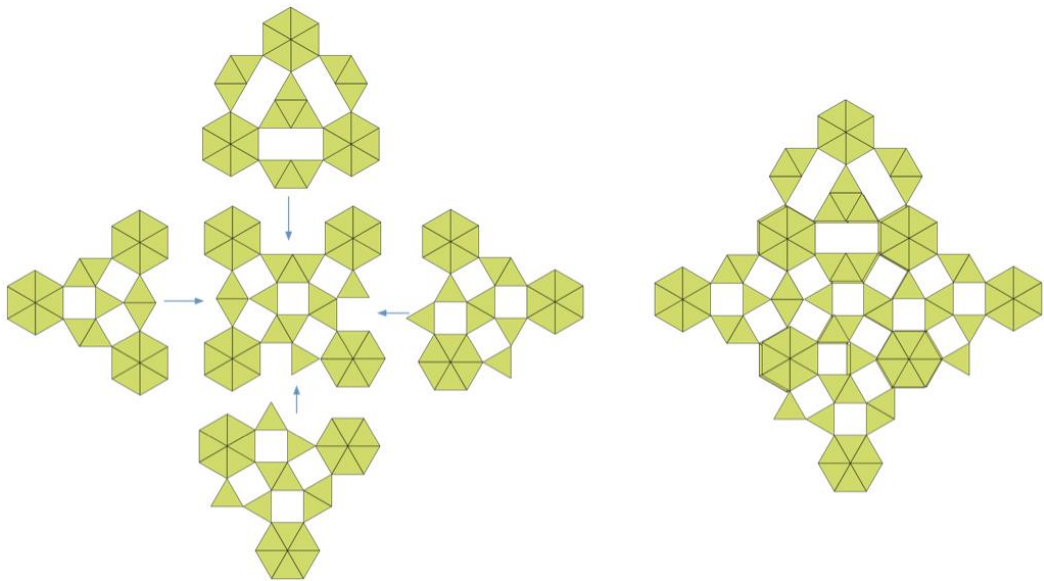

*Figure 7: Example of assembly of sets of triangles. An aperiodic network of triangles is constructed by applying the substitution rules. The vertices that connect two triangles are pivots that enable the empty spaces to take the shape of a rhombus of any angle.*

# Substitution rules for the aperiodic tiling made of a triangle and a rhombus

Let us know consider a rhombus of 60 degree angle. The aperiodic tiling made of triangles and 60-degree rhombi can be constructed with fixed substitution rules; these are given here below.

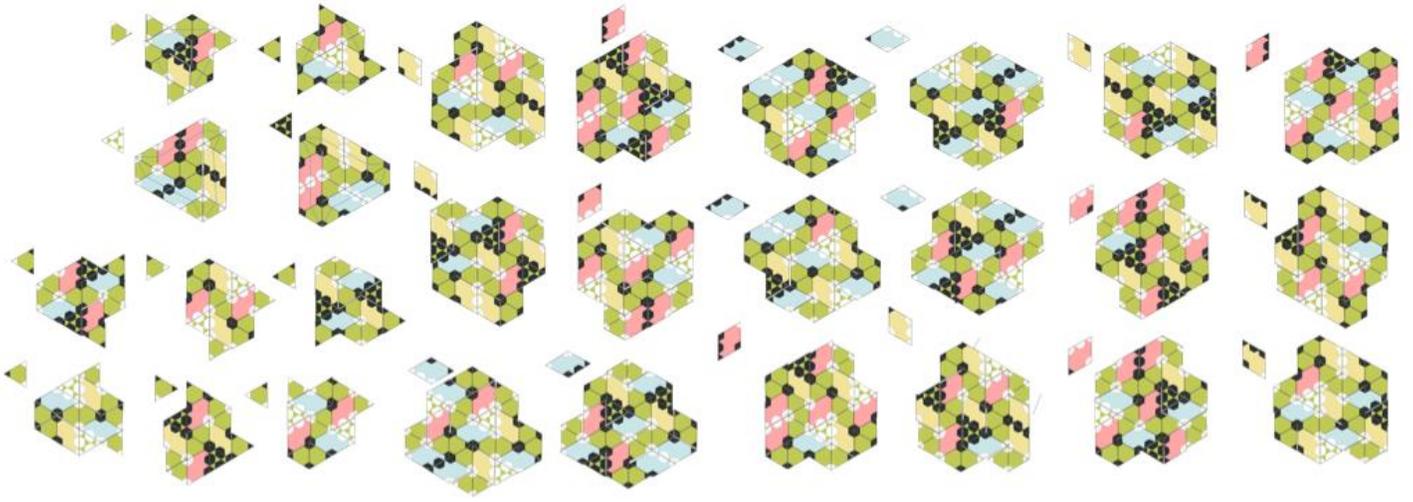

*Figure 8: Substitution rules to generate a tiling made of two tiles: a triangle and a rhombus of 60-degree angle.*

There are many rules compared to the square-triangle tiling but here the rotation of each tile is limited to 0 or 180 degrees, and no flipping is allowed. Here is an example of application of these tiles.

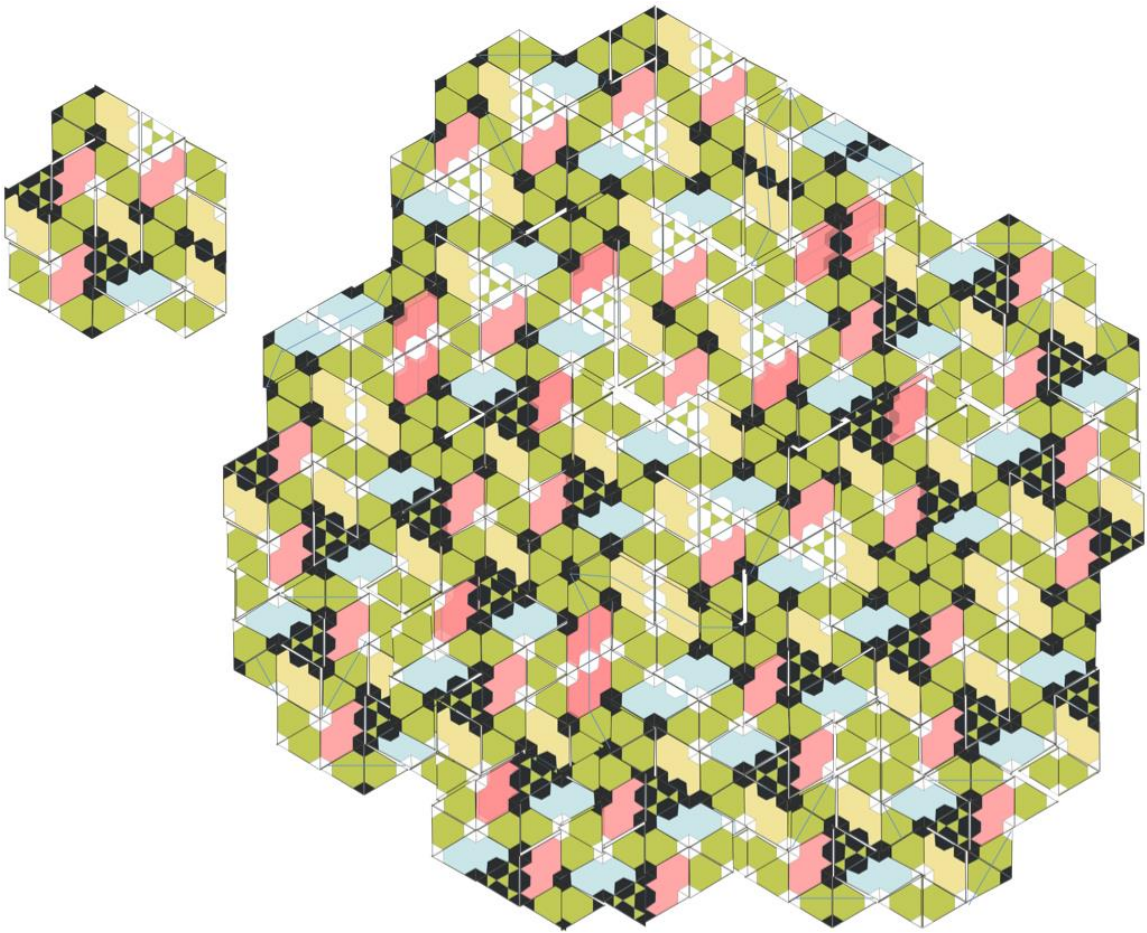

Its underlying structure obtained by connecting the centers of the hexagonal shapes (ie. made of 6 triangles) is shown here below. The structure is composed of triangles and rhombi. The size of these underlying rhombi are studied in the next section.

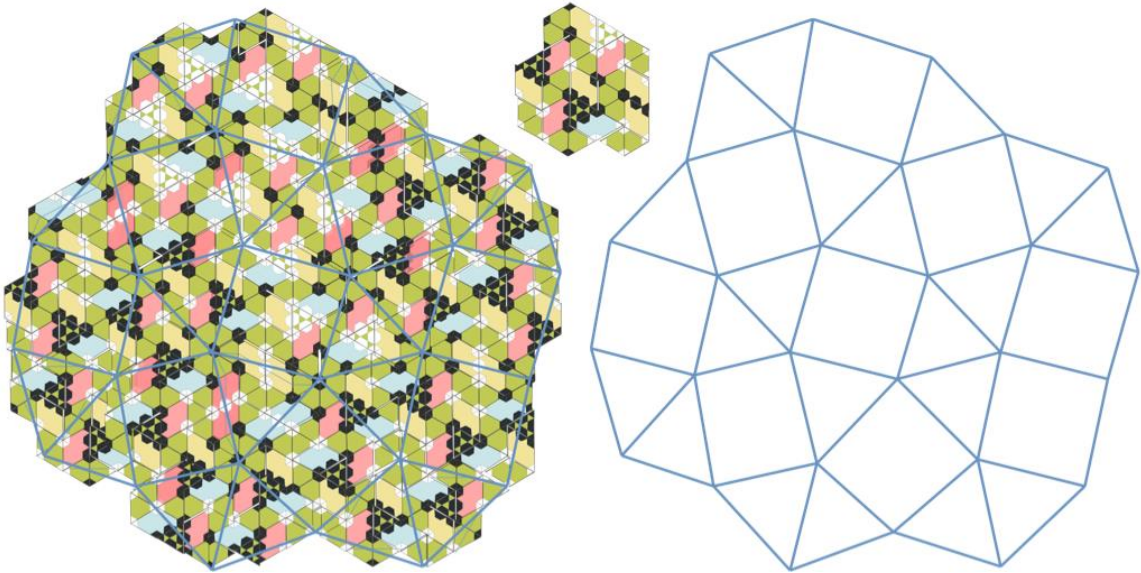

*Figure 10: The underlying structure of the aperiodic tiling made of two tiles: a triangle and a 60-degree rhombus.*

By allowing rotation of the tiles, we can define less substitution rules. Now tiles can rotate by 60 degrees. The number of rules is reduced to 12. They are shown here below.

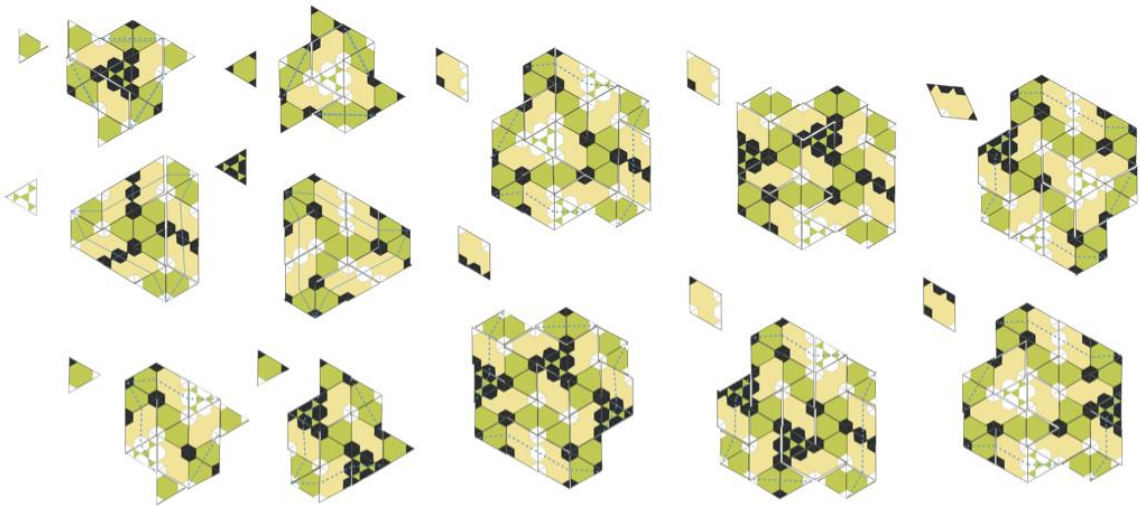

*Figure 11: Substitution rules with 60-degree rotation for the aperiodic tiling made of a triangle and a rhombus.*

In the above, each group of tiles has a negative counterpart that can be obtained by means of a flip and by changing the black into white and vice-versa. This allows us to further reduce the number of rules. In order to make the distinction between a tile and its flipped version, arrows are being used; All black corners have an arrow in one direction that is opposite to the one used for the white corners. Here are the 6 substitution rules that can be used to build an aperiodic tiling with a triangle and a 60-degree rhombus. Note that, compared to the previous ones shown above, these were further optimized to remove any overlap of tiles when assembling them.

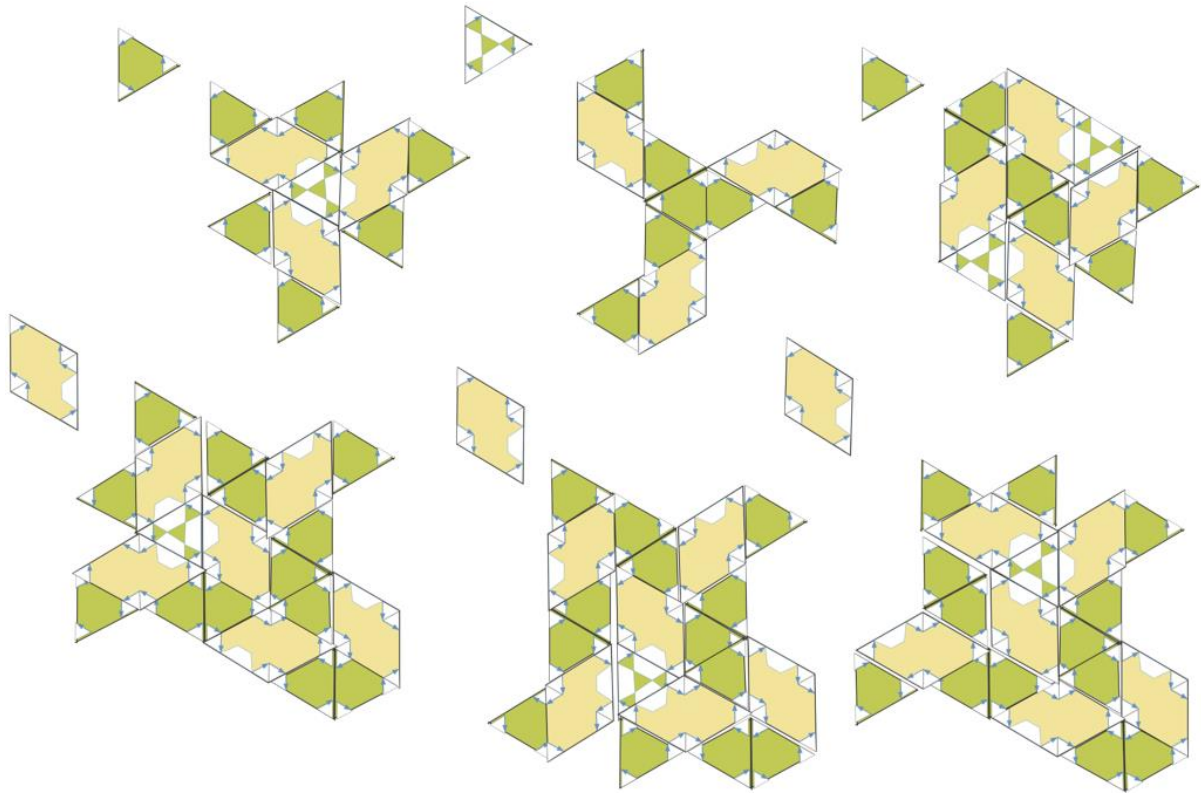

*Figure 12: Reduced set of substitution rules allowing rotation of 60 degrees and flipping.*

As explained before, the rhombus can be of any angle. For any given angle, those 6 rules can be adapted accordingly. Here below are the rules when the angle is 36 degrees.

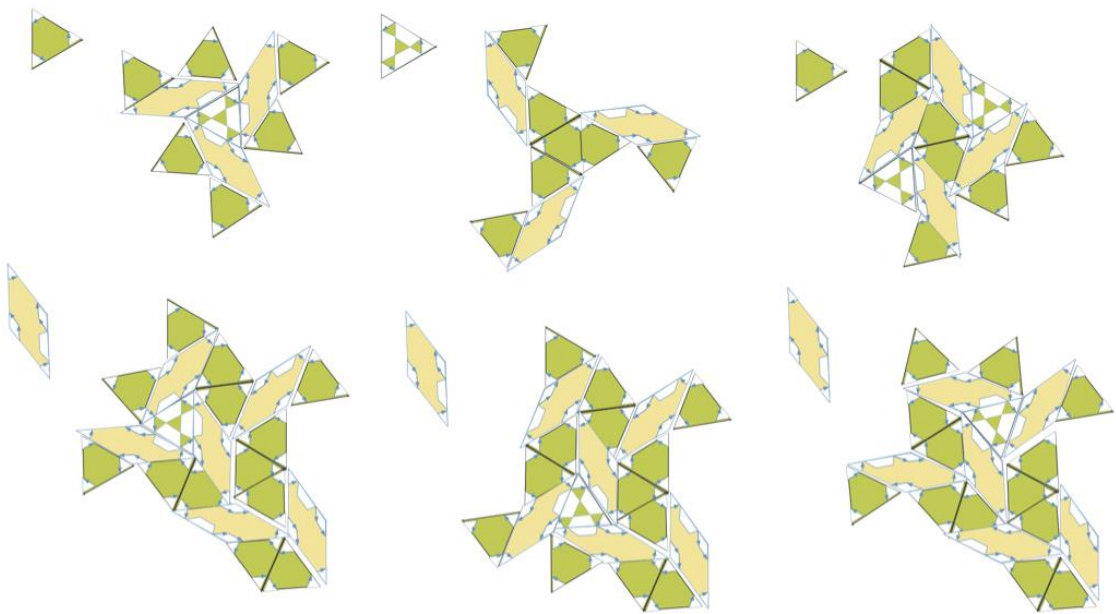

*Figure 13: Substitution rules for an aperiodic tiling made of two tiles: a triangle and a 36-degree rhombus.*

When using these rules, the following aperiodic tilings can be obtained: for a 60-degree rhombus and a 36-degree rhombus.

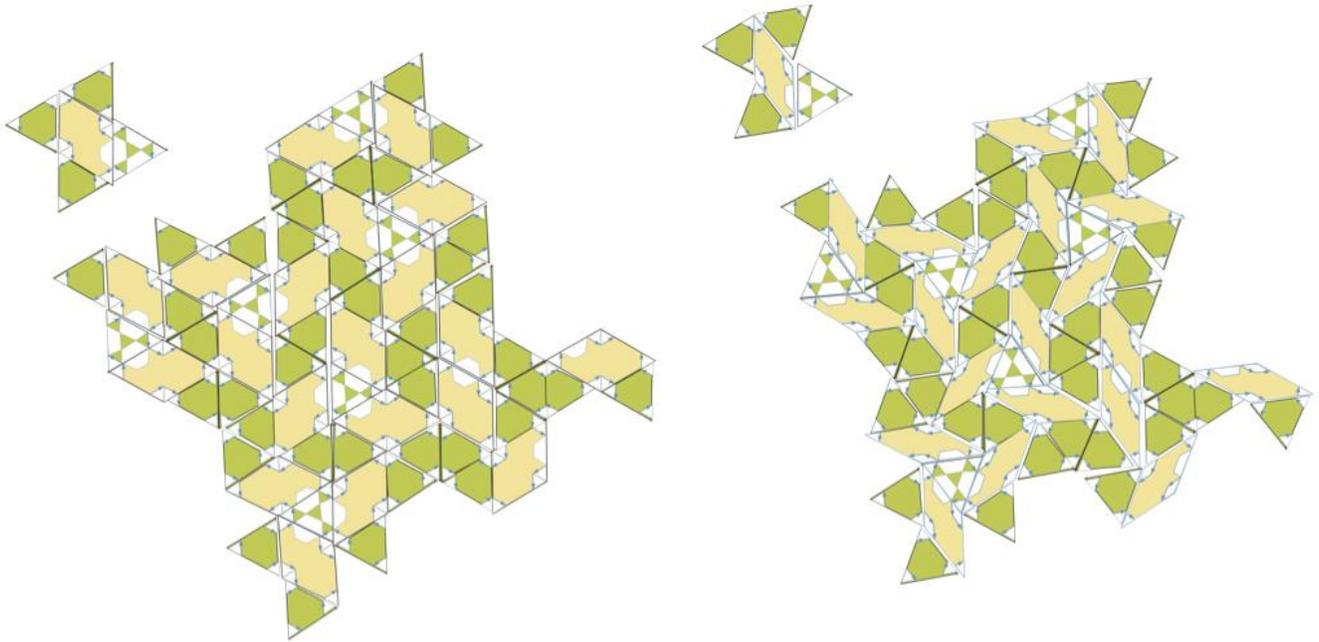

*Figure 14: Example of use of the substitution rules when the angle of the rhombus is 60-degree vs 36-degree.*

We will now further analyze these angles by taking another example, one where the rhombus has a 10-degree angle. Let us compare the effect of this change on two of the substitution rules that were used for the 60-degree rhombus. At the right of the figure, we show what happens to the new rules. Some of the triangles didn't change direction (the ones in green), others rotated by 10 degrees (the ones in red).

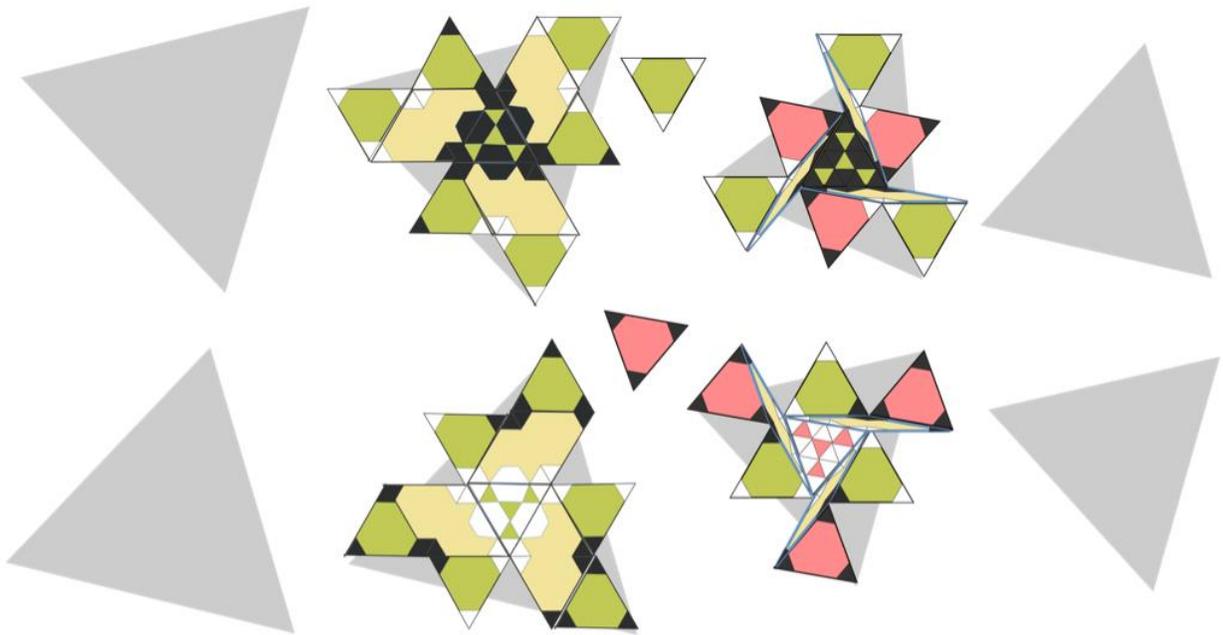

*Figure 15: Study of the angles when the rhombus is of 60 degree vs 10 degree.*

Based on this new coloring scheme, we propose the following 12 substitution rules. Here 60-degree rotations are needed but not flipping.

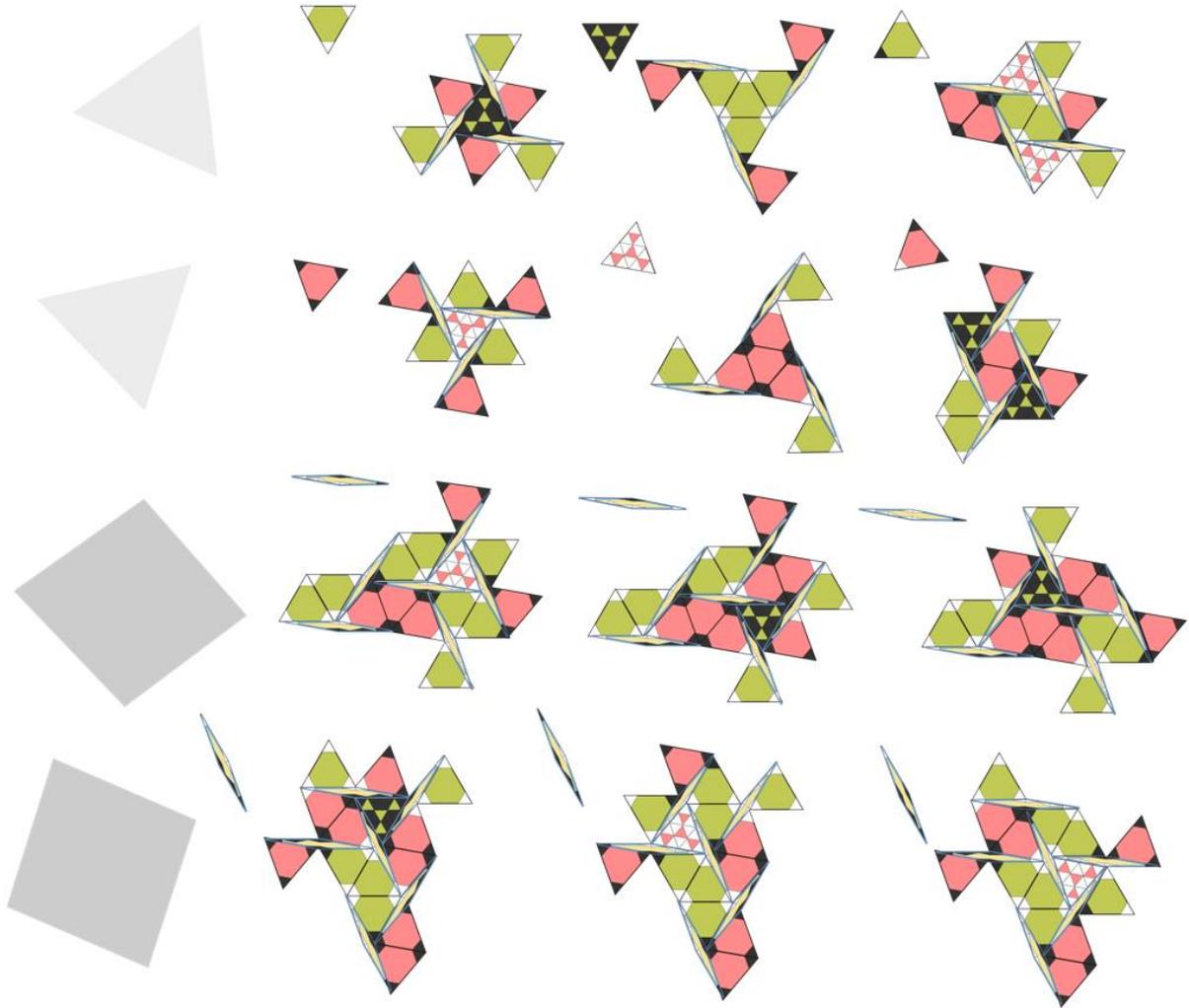

*Figure 16: Substitution rules for an aperiodic tiling made of two tiles: a triangle and a 10-degree rhombus. The use of these rules allows a rotation of 60 degrees but no flipping. All groups of tiles on the same line have the same underlying shape provided on the left (in grey): either an equilateral triangle or a rhombus of different angle.*

Here is an example of an aperiodic tiling produced with these rules.

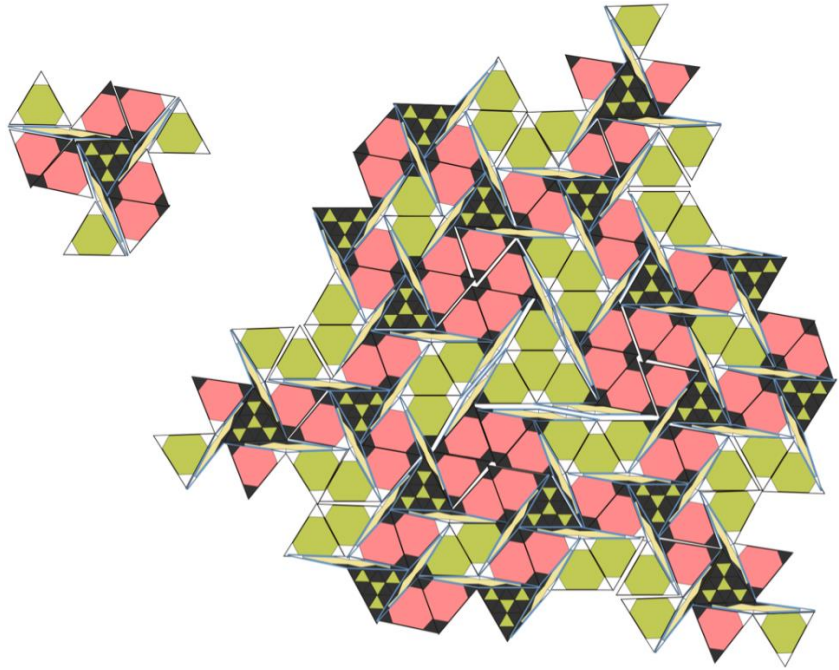

*Figure 17: Aperiodic tiling made of two tiles: a triangle and a rhombus of 10-degree angle.*

When connecting the center of the hexagons, the following underlying structure appears.

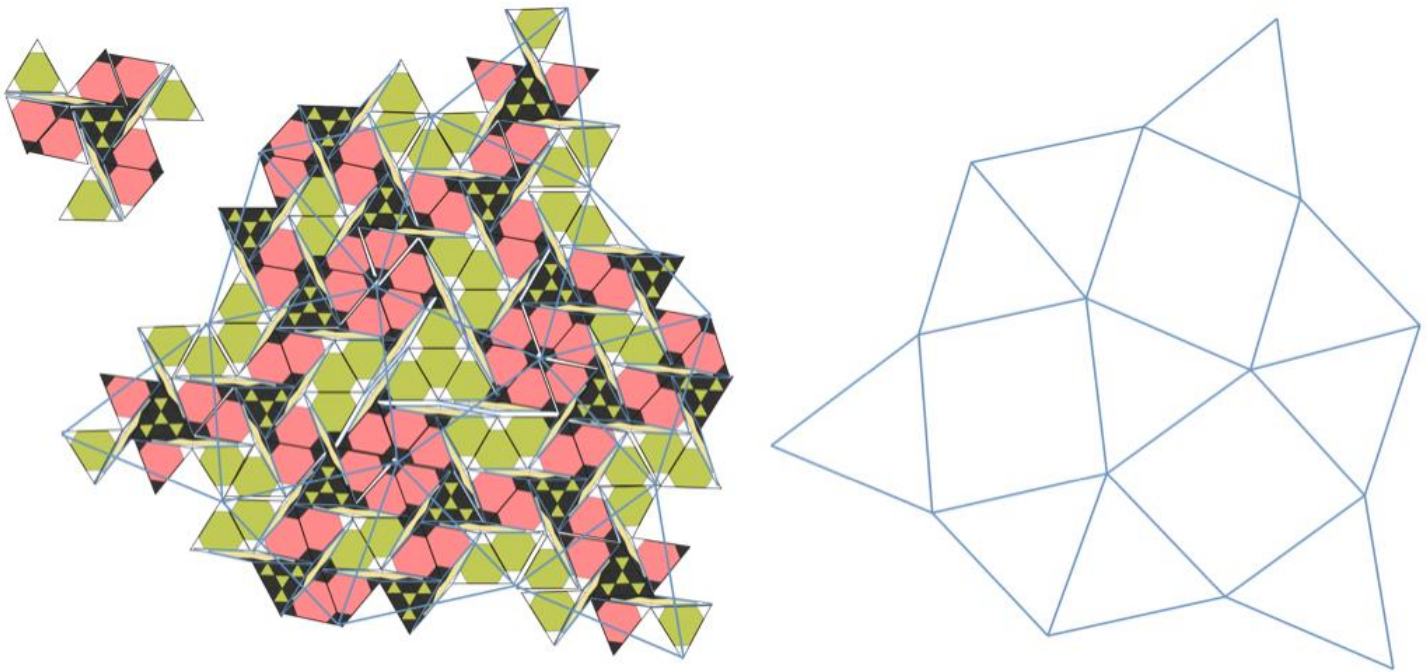

*Figure 18: Underlying structure when the rhombus is of 10-degree angle.*

Such aperiodic tiling can be created with a rhombus of any angle. The substitution rules can simply be adapted, or parameterized rules can be used, where the parameter is the angle of the rhombus. For any given angle, there will be green triangles whose direction is independent on that parameter, and red triangles that will rotate by the value of that angle.

Here is a summary on the different substitution rules provided in this paper for triangle-rhombus aperiodic tiling.

| Number of substitution rules | Possible rotations | Flipping |
|---|---|---|
| 28 | 180 degrees | Not allowed |
| 12 | 60 degrees | Not allowed |
| 6 | 60 degrees | Allowed |

## The limit case when the angle of the rhombus is near zero

Now, let us briefly discuss the limit case when the angle is near zero. It can be so thin that the rhombi become invisible. By using the same colors and motives on the tiles, we then obtain a nice aperiodic pattern made of one single tile that is visible: the triangle. If we remove the patterns on the tiles, we will simply see a periodic tiling of triangles.

Here below is the case when the rhombi have a 10-degree angle. The rhombi are like thick edges; they act like a glue to keep the triangular parts of the tiling together. The larger the angle of the rhombus is, the more apparent the glue becomes and the aperiodicity as well.

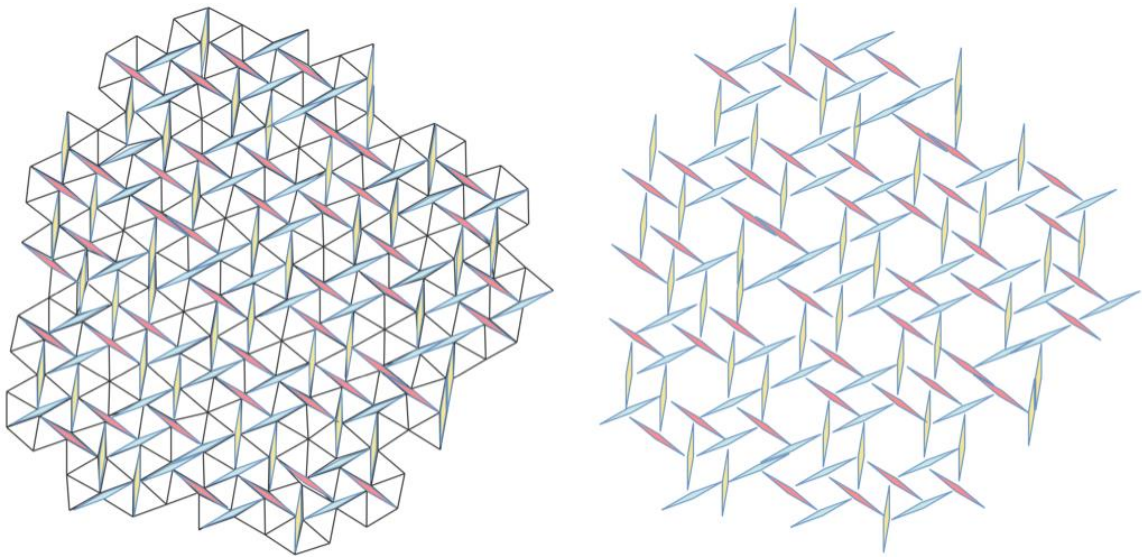

*Figure 19: Thin rhombi are like thick lines; they act like a glue to keep the triangular parts together. When the angle of the rhombus increases, the glue is more apparent and the aperiodicity as well.*

## Conclusion

Aperiodic tiling are interesting to study not only for their artistic value but also for their industrial applications. In this paper, we presented a new family of aperiodic tilings made of two tiles: triangles and rhombus. Their number is infinite as the rhombus can take any angle in the range 0 to 180 degrees. There are special cases. When the angle is 90 degree, the tiling is the aperiodic square-triangle one used in the study of dodecagonal quasicrystals. When the angle is either 0 or 180, the tiling becomes periodic. This set of triangle-rhombus tilings unifies aperiodic and periodic ones.

Another way of looking at the set is by considering the triangle-rhombus tiling to be one tiling of variable geometry. This by itself can lead to interesting new applications. If one modifies any rhombus of the tiling at a given location, all other rhombi will need to follow the same modification.

We showed that the triangle-rhombus aperiodic tilings have an underlying dodecagonal structure. Their use for the study of dodecagonal quasicrystals now needs further investigation.

## Special thanks

I would like to thank all the members of my family for being so supportive with my work. I am also thankful to Pierre Gradit, a researcher with whom I used to work some 30 years ago; this was a great opportunity to reconnect. I am also grateful to Professor U. Grimm who endorsed me on ArXiv. Finally, thanks to Marianne Impéror and her colleagues of U. Paris-Saclay for making me realize the potential use of aperiodic tilings of variable geometry in the research of aperiodic soft matters.
## References

Baake, M., & Grimm, U. (2013). *Aperiodic Order - Vol. 1 A Mathematical Invitation.* Cambridge University Press.

Baake, M., & Grimm, U. (2017). *Aperiodic Order - Vol. 2 Crystallography and Almost Periodicity.* Cambridge University Press.

Baake, M., Klitzing, R., & Schlottmann, M. (1992). Fractally shaped acceptance domains of quasiperiodic square-triangle tilings with dedecagonal symmetry. *Physica A: Statistical Mechanics and its Applications 191(1-4):554-558*, 554-558.

Grünbaum, B., & Shephard, G. C. (1986). *Tilings and Patterns.* New York: Freeman. ISBN 0-7167-1193-1.

Impéror-Clerc, M., Jagannathan, A., Kalugin, P., & Sadoc, J.-F. (2021). Square-triangle tilings : An infinite playground for softmatter. *Soft Matter,*, DOI: 10.1039/D1SM01242H.

Jian, L., & Yang, P. (2020). Antenna optimized array based on Schlottmann aperiodic tiling. *Asia-Pacific Conference on Microwave.* IEEE Xplore.

Oppenheim, J., Skorupskii, G., & Dincă, M. (2020). Aperiodic metal–organic frameworks. *Chemical Science*.

Penrose, R. (1979). Pentaplexity A Class of Non-Periodic Tilings of the Plane. *The Mathematical Intelligencer 2 (https://doi.org/10.1007/BF03024384)*, 32–37 .

Sadoc, J.-F., & Mosseri, R. (n.d.). Square-triangle tilings: a template for soft-matter quasicrystals. *http://www.xray.cz/aperiodic/files/26.htm*.

Shechtman, D. (n.d.). Retrieved from https://en.wikipedia.org/wiki/Dan_Shechtman

Steurer, W. (2004). Twenty years of structure research on quasicrystals. Part I. Pentagonal, octagonal, decagonal and dodecagonal quasicrystals". *Zeitschrift für Kristallographie - Crystalline Materials*, 391-446 https://doi.org/10.

Varn, D., & Crutchfield, J. (2016). What did Erwin mean? The physics of information from the materials genomics of aperiodic crystals and water to molecular information catalysts and life. *Philosophical Transactions of the Royal Society A*.

Weisstein, E. W. (n.d.). *Aperiodic Tiling*. Retrieved from MathWorld--A Wolfram Web Resource. : https://mathworld.wolfram.com/AperiodicTiling.html

Wikipedia. (n.d.). *Penrose Tiling*. Retrieved from https://en.wikipedia.org/wiki/Penrose_tiling
13